\documentclass[12pt]{amsart}
 
\usepackage{epsfig}
\usepackage{graphics}
\usepackage{latexsym}
\usepackage{psfrag}
\usepackage{amssymb}
\usepackage{amsmath}
\usepackage{amsthm} 

\newcommand{\cc}{{\mathcal{C}}}

\renewcommand{\setminus}{{\smallsetminus}}
\newcommand{\NN}{{\mathbb{N}}}
\newcommand{\BB}{{\mathbb{B}}}

\newcommand{\HH}{{\mathbb{H}}}

\newcommand{\nin}{{\notin}}

\newcommand{\cover}[1]{{\widetilde{#1}}}

\newcommand{\lift}[1]{{\widetilde{#1}}}
\newcommand{\closure}[1]{{\overline{#1}}}

\newcommand{\ML}{{\mathcal{ML}}} % Ditto for laminations
\newcommand{\PML}{{\mathcal{PML}}}

\newcommand{\homeo}{{\medspace \cong \medspace}} 
 
\newcommand{\neigh}{{\eta}} % Open regular neighborhood
\newcommand{\bdy}{{\partial}} % Boundary

\newcommand{\RRPP}{{\mathbb{RP}}}

\newcommand{\genus}{{\rm{genus}}}

\theoremstyle{plain}
\newtheorem{theorem}{Theorem}[section]

\newtheorem{lemma}[theorem]{Lemma}

\theoremstyle{definition}
\newtheorem*{define}{Definition}

\newtheorem{proposition}[theorem]{Proposition}
\newtheorem{remark}[theorem]{Remark}
\newtheorem*{example}{Example}

\newsavebox{\savepar}

\begin{document}

\title{Strongly Irreducible Surface Automorphisms}
\author{Saul Schleimer}
\address{\hskip-\parindent
        Saul Schleimer\\
        Department of Mathematics, UIC\\
        851 South Morgan Street\\
        Chicago, Illinois 60607}
\email{saul@math.uic.edu}

\date{\today}

\begin{abstract}
A surface automorphism is {\em strongly irreducible} if every
essential simple closed curve in the surface has nontrivial geometric
intersection with its image.  We show that a three-manifold admits
only finitely many inequivalent surface bundle structures with
strongly irreducible monodromy.
\end{abstract}
\maketitle

\section{Introduction}

A surface automorphism $h:F \rightarrow F$ is {\em strongly
irreducible} if every essential simple closed curve $\gamma \subset F$
has nontrivial geometric intersection with its image, $h(\gamma)$.
This paper shows that a three-manifold admits only finitely many
inequivalent surface bundle structures with strongly irreducible
monodromy.  This imposes a serious restriction; for example, any
three-manifold which fibres over the circle and has $b_2(M) \geq 2$
admits infinitely many distinct surface bundle structures.

The main step is an elementary proof that all {\em weakly
acylindrical} surfaces inside of an irreducible triangulated manifold
are isotopic to fundamental normal surfaces.  As weakly acylindrical
surfaces are a larger class than the acylindrical surfaces this
strengthens a result of Hass~\cite{Hass95}; an irreducible
three-manifold contains only finitely many acylindrical surfaces.

Section~\ref{DefinitionsExamples} gives necessary topological
definitions, examples of strongly irreducible automorphisms, and
precise statements of the theorems.  The required tools of normal
surface theory are presented in Section~\ref{NormalSurfaces}.
Section~\ref{WeaklyAcylSurfaces} defines weakly acylindrical and proves that
every such surface is isotopic to a fundamental surface.  In the spirit
of the Georgia Topology Conference the paper ends by listing several
open questions.

Many of the ideas and terminology discussed come from the study of
Heegaard splittings as in~\cite{CassonGordon87} and in my thesis
\cite{Schleimer01}.  This paper, in particular
Theorem~\ref{WeaklyAcylImpliesFundamental}, owes an obvious debt
to~\cite{JacoOertel84} by Jaco and Oertel.  I thank Ian Agol for
simplifying my original proof of
Proposition~\ref{PeriodicStronglyIrreducible} and Dave Bachman for
critiquing an early version of this paper.

\section{Definitions and Examples}
\label{DefinitionsExamples}

This section lays out the necessary definitions, states the main
theorems precisely, and gives examples of families of strongly
irreducible surface automorphisms.

Let $F$ be a closed, orientable, $\genus(F) > 1$ surface.  Let $h:F
\rightarrow F$ be an automorphism of $F$.  If $\gamma_0$ and
$\gamma_1$ are simple closed curves in $F$ then the {\em geometric
intersection number}, $i(\gamma_0, \gamma_1)$, is the minimum of
$|\gamma'_0 \cap \gamma'_1|$ taken over all $\gamma'_i$ isotopic to
$\gamma_i$.

\begin{define}
The map $h$ is {\em strongly irreducible} if $i(h(\gamma), \gamma) >
0$ for every essential simple closed curve $\gamma \subset F$.
\end{define}

If $h$ is not strongly irreducible then $h$ is {\em weakly reducible}.

\begin{remark}
A {\em reducible} surface automorphism is one which admits an
invariant set of disjoint essential simple closed curves, up to
isotopy.  Thus reducible maps are also weakly reducible.
\end{remark}

\begin{remark}
There exist weakly reducible pseudo-Anosov maps with arbitrarily high
stretch factor --- the main ideas required for the construction may be
found in~\cite{LongMorton86}.
\end{remark}

\begin{example}
Let $P$ be a regular $4g$-gon in the hyperbolic plane with angle at
the vertices equal to $2\pi/4g$.  Here $g$ is assumed to be two or
larger.  Glue opposite sides of $P$ by isometries to obtain $F$, an
orientable surface of genus $g$.  Let $h:P \rightarrow P$ be a
counter-clockwise rotation of $P$ about its center, $O$, through an
angle of $2\pi/4g$.  Let $h'$ be the induced isometry of $F$.

\begin{proposition}
\label{PeriodicStronglyIrreducible}
The periodic map $h'$ is strongly irreducible.
\end{proposition}

All of the hyperbolic trigonometry needed in the proof may be found in
Chapter~7 of~\cite{Beardon83}.

\begin{proof}
Suppose that $R$ is the distance between $O$ and $V$, where $V$ is a
vertex of $P$.  Then $\cosh(R) = \cot^2(2\pi/8g)$.  Both $O$ and $V$
give fixed points, $x_O, x_V \in F$, of $h'$.  Suppose that $\gamma
\subset F$ is a simple closed geodesic.  There is a point of $\gamma$
which lies within $R/2$ of either $x_O$ or $x_V$.  This last holds
because the points of $F$ which are not this close to one of $x_O$ or
$x_V$ form a disjoint union of disks.

Suppose there is a point of $\gamma$ within distance $R/2$ of the
point $x_O$.  (The other case is similar.)  Pick $\lift{\gamma}
\subset \HH^2$, a lift of $\gamma$, which lies within $R/2$ of $O$.
Let $L$ be the distance between $\lift{\gamma}$ and $O$.  Let $\theta$
be the visual angle which $\lift{\gamma}$ occupies, as viewed from
$O$.  (When $\lift{\gamma}$ meets $x_O$ take $\theta = \pi$.)  Then $L
\leq R/2$ and $\cosh(L) = 1/\sin(\theta/2).$ Applying a
``double-angle'' formula to $\cosh(R)$ yields
$$\cosh(R/2) = 1/\sqrt{2}\sin(2\pi/8g).$$
Hyperbolic cosine is an increasing function on the positive reals so
$$1/\sin(\theta/2) \leq 1/\sqrt{2}\sin(2\pi/8g).$$ As sine is
increasing in the interval $[0, \pi/2]$ we deduce that $\theta >
2\pi/4g$.  Thus the visual angle of $\lift{\gamma}$ is greater than
$2\pi/4g$ and $h(\lift{\gamma}) \cap \lift{\gamma} \neq \emptyset$.
It follows that $h'(\gamma) \cap \gamma \neq \emptyset$.  This gives
the desired conclusion as the intersection and geometric intersection
numbers agree for geodesics.
\end{proof}
\end{example}

\begin{proposition}
\label{PAStronglyIrreducible}
If $h$ is a pseudo-Anosov map then there is an $n \in \NN$ such that
$h^n$ is strongly irreducible.
\end{proposition}

Here we only sketch a proof; the numbers in parentheses refer to
theorems in Kapovich's book~\cite{Kapovich01} which we take as our
reference for $\PML(F)$, the space of projectively measured
laminations of $F$.  Recall that geometric intersection extends to a
continuous function on $\ML \times \ML$ (11.26).

Let $\lambda^\pm$ be the stable and unstable laminations for $h$.  Let
$U, V$ be small neighborhoods of $\lambda^\pm$ (respectively) in
$\PML(F)$.  Choose $U$ and $V$ so that all $x \in U$, $y \in V$ have
$i(x,y) > 0$.  This is possible because the geometric intersection
between $\lambda^+$ and $\lambda^-$ is non-zero (11.49).  In
particular, $U \cap V = \emptyset$.

There is an $m \in \NN$ such that if $x \in \PML(F) \setminus V$ then
$h^m(x) \in U$ (11.47).  Now, $n = m + 1$ gives the desired
conclusion.  To see this, pick an essential simple closed curve $y
\subset F$.  Denote the corresponding element of $\PML(F)$ again by
$y$.  There is a integer $k$ such that $h^k(y) \in V$ but $h^{k+1}(y)
\nin V$.  (Again,~11.47.)  So $i(y, h^n(y)) = i(h^k(y), h^{k+1+m}(y))
> 0$.  This completes the proof sketch.

\begin{remark}
\label{MasurMinsky}
The notions irreducible and strongly irreducible may be generalized by
introducing the {\em curve complex} of $F$, $\cc(F)$, and defining the
translation distance, $\tau(h)$, of $h$'s action on $\cc(F)$.  The
deep results of~\cite{MasurMinsky99} give a positive integer $n(g)$
depending only on $g = \genus(F)$ such that: If $h$ is pseudo-Anosov
then $\tau(h^{n(g)}) \geq 2$ and so $h^{n(g)}$ is strongly
irreducible.  This greatly improves upon
Proposition~\ref{PAStronglyIrreducible}.
\end{remark}

We now turn from examples to the main objects of interest: surface
bundles over the circle.

\begin{define}
If $h:F \rightarrow F$ is a surface automorphism then let $M_h$ be the
{\em mapping torus} of $h$.  So
$$M_h \homeo F \times I / {(x,1) \sim (h(x),0)}.$$
\end{define}
The mapping torus admits a natural map to the circle, $\pi_h : M_h
\rightarrow S^1$.  The embedded surfaces $F \times \{t\} \subset M$
are called the {\em fibres} of this map while $h$ is the {\em
monodromy} and $g(F) = \genus(F) \geq 2$ is the {\em genus} of the
bundle.  As an example, if $h$ is a strongly irreducible periodic
automorphism, then $M_h$ is an atoroidal Seifert fibred space.

Fix $M$, a closed orientable three-manifold.

\begin{define}
A {\em surface bundle structure} on $M$ is a pair $(h, \phi)$, where
$h$ is a surface automorphism and $\phi$ is a homeomorphism from $M_h$
to $M$.
\end{define}

Let $(h,\phi)$ and $(h',\phi')$ be two surface bundle structures on
$M$.  Suppose that $\psi : M_h \rightarrow M_{h'}$ is a homeomorphism
such that $\pi_h = \pi_{h'} \circ \psi$ and $\phi$ is isotopic to
$\phi' \circ \psi$.  Then the two bundle structures on $M$ are {\em
equivalent}.  Thus $(h, \phi)$ is equivalent to $(h', \phi')$ if and
only if $h'$ is conjugate to $h$ and the fibres of the two bundle
structures are isotopic in $M$.

Here is a precise statement of the theorem alluded to in the
introduction. 

\vspace{2mm}
\noindent
{\bf Theorem~\ref{A}.} 
{\em Suppose that $M$ is a closed, orientable three-manifold.  Then $M$
admits only finitely many inequivalent surface bundle structures with
strongly irreducible monodromy.}
\vspace{1.5mm}

Closely related is our:

\vspace{2mm}
\noindent
{\bf Theorem~\ref{B}.} 
{\em Suppose that $M$ is a closed, orientable three-manifold.  There is a
positive real number $c(M)$ such that if $(h,\phi)$ is a surface bundle
structure on $M$ with genus $g > 1$ then $h^i$ is weakly reducible
for all integers $i$ where $1 \leq i \leq c(M) \cdot (2 g - 2)$.}
\vspace{1.5mm}

If $c(M) \cdot (2g - 2) < 1$ then the theorem is vacuous as no
such $i$ exists. 

\begin{remark}
If $M$ is atoroidal then $M$ has only finitely many inequivalent
surface bundle structures in each genus.  (This simply because there
are only finitely many incompressible surfaces, up to isotopy, of each
genus.  See~\cite{JacoOertel84}.)  Thus, in the atoroidal case,
Theorem~\ref{B} implies Theorem~\ref{A}.
\end{remark}

\begin{remark}
\label{LowerBoundOnN(g)}
Suppose $M$ is a closed, atoroidal three-manifold admitting infinitely
many inequivalent surface bundle structures.  Then this manifold and
Theorem~\ref{B} provide another proof that the constant $n(g)$ of
Remark~\ref{MasurMinsky} must tend to infinity as $g$ does.
\end{remark}

\begin{remark}
\label{Hartshorn}
Again, consider the translation distance $\tau(h)$ where $h$ is an
automorphism of the closed, orientable, genus at least two, surface
$F$.  One obtains an analogue of a theorem of
Hartshorn~\cite{Hartshorn99} regarding Heegaard splittings: If
$(h,\phi)$ is a surface bundle structure on $M$ and $G$ is a two-sided
incompressible surface in $M$ then either $G$ is isotopic to a fibre,
$G$ is a torus (and $h$ is reducible), or $\tau(h) \leq - \chi(G)$.
 
As a corollary, deduce that if $h$ has translation distance greater
than $- \chi(F)$, where $F$ is the fibre, then $F$ is the unique
minimal genus incompressible surface in $M$, up to isotopy.  Weak
conclusions about the shape of the Thurston norm ball and on the
structure of $M$'s symmetry group follow.
\end{remark}

\section{Normal surfaces}
\label{NormalSurfaces}

\begin{figure}
$$\begin{array}{cc}
\epsfig{file=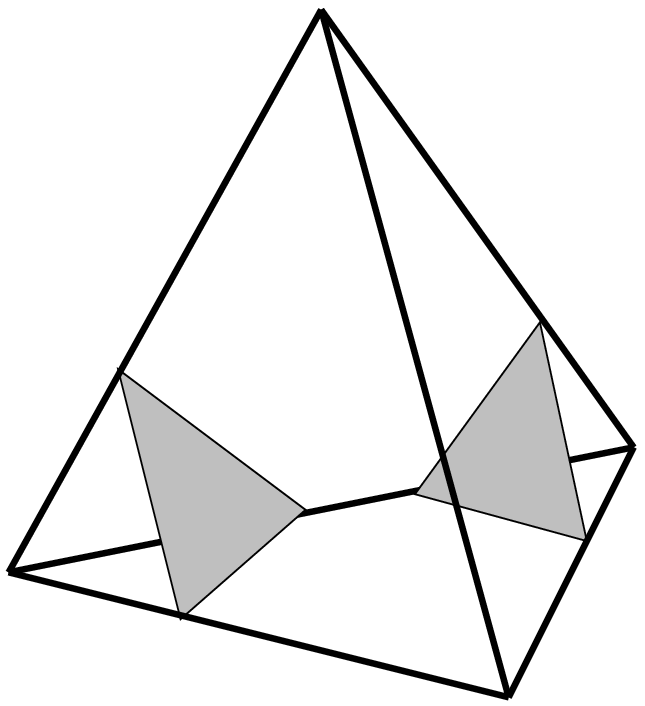, height=3.5cm} &
\epsfig{file=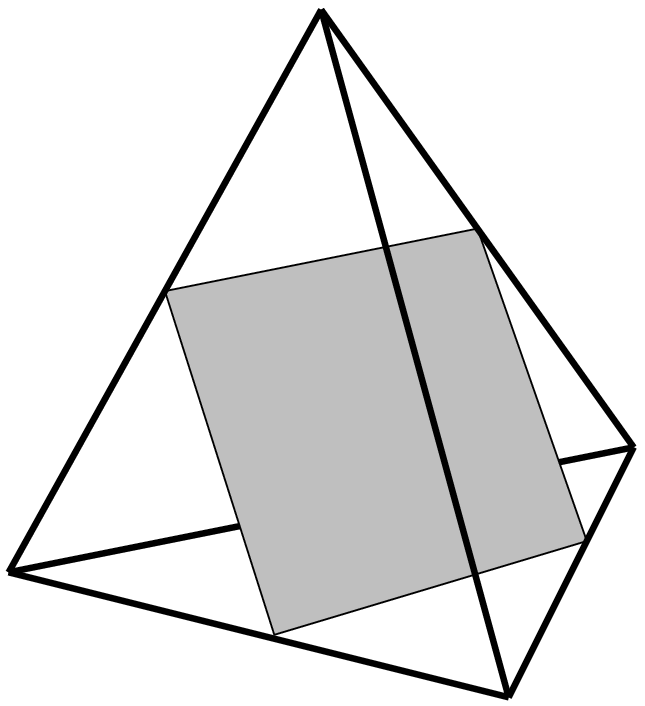, height=3.5cm}
\end{array}$$
\caption{Normal triangles and normal quad}
\label{NormalDisks}
\end{figure}

This section presents the required bare minimum of normal surface
theory.  For a more complete treatment consult~\cite{JacoOertel84}
or~\cite{JacoTollefson95}.

Fix a closed orientable three-manifold $M$ and choose $T$, a
triangulation of $M$.  Suppose $F \subset M$ is a closed embedded
surface.  The {\em weight} of $F$, $w(F)$, is the number of
intersections between $F$ and the one-skeleton of $T$.  The surface
$F$ is {\em normal} if $F$ intersects every tetrahedron of $T$ in a
disjoint collection of {\em normal triangles and quadrilaterals}.  See
Figure~\ref{NormalDisks}.  In each tetrahedron there are four types of
normal triangle and three types of normal quad.

\begin{lemma}{\rm (Haken~\cite{Haken68})}
\label{IncompressibleImpliesNormal}
Suppose that $(M, T)$ is closed and irreducible.  If $F \subset M$ is
embedded and incompressible then $F$ is isotopic to a normal surface
$F'$ with $w(F') \leq w(F)$.
\end{lemma}

A {\em normal isotopy} of $M$ fixes every simplex of $T$ setwise.  Two
normal surfaces, $G$ and $H$, are {\em compatible} if in each
tetrahedron $G$ and $H$ have the same types of quad (or one or both
have no quads.)  In this situation we form the {\em Haken sum} $F = G +
H$ as follows:

Normally isotope $H$ to make $G$ transverse to $H$, to make $\Gamma =
G \cap H$ transverse to the skeleta of $T$, and to minimize the number
of curves in $\Gamma$.  The components of $\Gamma$ are the {\em
exchange curves}.  For every such $\gamma \subset \Gamma$ let
$R(\gamma)$ be the closure (taken in $M$) of $\neigh_M(\gamma)$, an
open regular neighborhood of $\gamma$.  The set $R(\gamma)$ is a solid
torus containing $\gamma$ as a core curve.

\begin{figure}
\psfrag{Ar}{$A_r$}
\psfrag{Ai}{$A_i$}
$$\begin{array}{c}
\epsfig{file=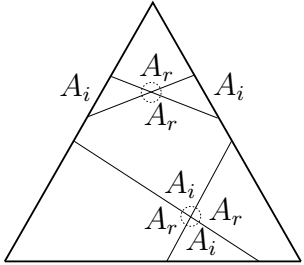, height=3.5cm}
\end{array}$$
\caption{Regular and irregular annuli}
\label{SmallNeighborhood}
\end{figure}

\begin{figure}
\psfrag{F}{$G$}
\psfrag{G}{$H$}
\psfrag{Arg}{$A_r(\gamma)$}
\psfrag{gamma}{$\gamma$}
\psfrag{regularexchange}{Cut-and-paste}
\psfrag{exchangeband}{$C(\gamma)$}
$$\begin{array}{c}
\epsfig{file=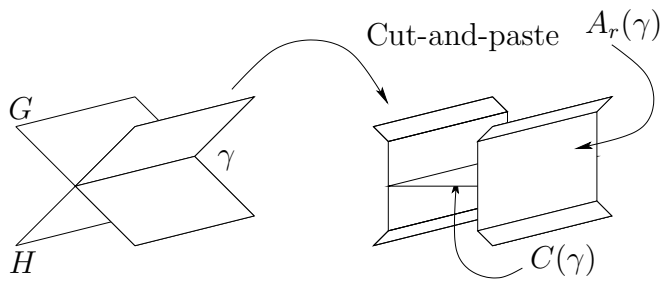, height=3.5cm}
\end{array}$$
\caption{The exchange band for $\gamma$.}
\label{ExchangeBand}
\end{figure}

Then $(\bdy R(\gamma)) \setminus (G \cup H)$ is a union of annuli.
Taking closures divide these into two sets, the {\em regular} and {\em
irregular} annuli, $A_r(\gamma)$ and $A_i(\gamma)$, as indicated by
Figure~\ref{SmallNeighborhood}.  Finally, as in
Figure~\ref{ExchangeBand}, form the surface

$$F = \left( (G \cup H) \setminus \bigcup_\Gamma \{\neigh_M(\gamma)\} \right)
             \cup \bigcup_\Gamma \{A_r(\gamma)\}.$$

Each connected component of $(G \cup H) \setminus \bigcup_\Gamma
\{\neigh_M(\gamma)\}$ is a {\em patch} of the sum $G + H$ while each
closed annulus ${A_r(\gamma)}$ is a {\em seam}.  Note that $F$ is again
a normal surface which, up to normal isotopy, does not depend on the
choices made in the above construction.

A normal surface is {\em fundamental} if it admits no such
decomposition.  A fundamental result due to Haken is:

\begin{lemma}{\rm (See~\cite{JacoTollefson95})}
\label{FiniteFundamental}
A closed orientable triangulated three-manifold $(M^3,T)$ contains
only finitely many fundamental normal surfaces, up to normal isotopy.
\end{lemma}

\begin{remark}
The proof of Lemma~\ref{FiniteFundamental} requires only that the
manifold $M$ admit a finite triangulation.  We restrict ourselves to
the closed and orientable case to avoid an unnecessarily technical
discussion about normal surfaces.
\end{remark}

Each ``cut-and-paste'' operation involved in the Haken sum $F = G + H$
may be recorded by an embedded {\em exchange band} $(C(\gamma), \bdy
C(\gamma)) \subset (N(\gamma), A_r(\gamma))$.  That is, the band is
embedded in $N(\gamma)$ with boundary inside of $A_r(\gamma)$.  See
Figure~\ref{ExchangeBand}.  Each $C(\gamma)$ is either an annulus or a
Mobius band.  The exchange bands record enough information to reverse
the sum.  Note also that each seam is a regular neighborhood (in $F$)
of a boundary component of some $C(\gamma)$.

A sum $F = G + H$ is {\em reduced} if $F$ cannot be realized as a sum
$G' + H'$ where $G'$ and $H'$ are again normal, isotopic to $G$ and
$H$ respectively, with $|G' \cap H'| < |G \cap H|$.  Note that the
isotopy between $G$ and $G'$ ($H$ and $H'$) need not be normal.
Lemma~\ref{NoDiskPatches} is a key technical result
for~\cite{JacoOertel84} and our
Theorem~\ref{WeaklyAcylImpliesFundamental}.

\begin{lemma}{\rm (Jaco and Oertel~\cite{JacoOertel84})}
\label{NoDiskPatches}
Suppose that $(M^3,T)$ is a closed, orientable, irreducible,
triangulated three-manifold.  Suppose that $F \subset M$ is an
incompressible normal surface which is least weight in its isotopy
class.  If the sum $F = G + H$ is reduced then no patch of $G + H$ is
a disk.
\end{lemma}

\section{Weakly acylindrical surfaces}
\label{WeaklyAcylSurfaces}

Here the {\em weakly acylindrical} surfaces are defined.  Equipped
with this definition and Theorem~\ref{WeaklyAcylImpliesFundamental} we
will prove Theorem~\ref{A}, the main goal of the paper.

Suppose now that $N$ is a compact, orientable three-manifold with
non-empty boundary.

\begin{define}
An embedded annulus $(A, \bdy A) \subset (N, \bdy N)$ is {\em
essential} in $N$ if $A$ is incompressible and
boundary-incompressible.
\end{define}

Now fix $M$, a closed orientable three-manifold.  Let $F \subset M$ be
a closed, embedded, incompressible, two-sided surface with genus at
least two.

\begin{define}
The surface $F$ is {\em cylindrical} if $N = M \setminus \neigh_M(F)$
admits an essential annulus.  If $N$ does not admit an essential
annulus then $F$ is {\em acylindrical}.
\end{define}

\begin{define}
The surface $F$ is {\em strongly cylindrical} if there exists an
embedded annulus $(A, \bdy A) \subset (M, F)$ such that $A \setminus
\neigh_M(F)$ is essential in $N = M \setminus \neigh_M(F)$.  On the
other hand, if no such annulus exists then $F$ is {\em weakly
acylindrical}.
\end{define}

\begin{remark}
\label{WeaklyAcyl}
If $F$ is acylindrical then $F$ is weakly acylindrical.  Note that if
$F$ is a fibre of a surface bundle then $F$ is never acylindrical, as
the complement of $F$ is homeomorphic to $F \times I$.  However, if
the monodromy is strongly irreducible then $F$ is weakly acylindrical.
For separating surfaces the notions acylindrical and weakly
acylindrical coincide.
\end{remark}

\begin{theorem}
\label{WeaklyAcylImpliesFundamental}
Suppose that $(M,T)$ is a closed, orientable, irreducible,
triangulated three-manifold.  Suppose $F \subset M$ is weakly
acylindrical.  Then $F$ is isotopic to a fundamental normal surface.
\end{theorem}

Before proceeding with the proof of
Theorem~\ref{WeaklyAcylImpliesFundamental} we remark that the normal
surface techniques used may easily be replaced by the methods of
branched surface theory.

\begin{figure}
\psfrag{F}{$F''$}
$$\begin{array}{c}
\epsfig{file=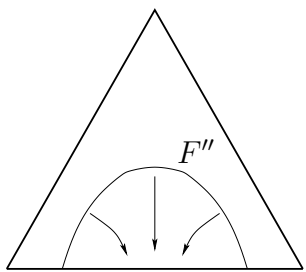, height = 3.5 cm}
\end{array}$$
\caption{The weight-reducing isotopy}
\label{BentArcWithIsotopy}
\end{figure}

\begin{proof}[Proof of Theorem~\ref{WeaklyAcylImpliesFundamental}]  
Fix $F \subset M$ as in the hypothesis.  Recall that weakly
acylindrical includes the property of being incompressible.
Applying Lemma~\ref{IncompressibleImpliesNormal} isotope $F$ to be
normal and least weight in its isotopy class.  Set $N = M \setminus
\neigh_M(F)$.  Recall also that $F$ is two-sided, and has genus two or
greater.

Suppose that $F$ is not fundamental.  Pick a reduced decomposition $F
= G + H$.  This sum admits some exchange band $A$.  Let $A_N = A \cap
N$.  If $A$ is a Mobius band then $A$'s {\em double}, defined below,
must be compressible or boundary-compressible in $N = M \setminus
\neigh_M(F)$.  If $A$ is an annulus then $A_N$ is itself compressible
or boundary-compressible in $N$.  The proof deals with each of these
possibilities in turn, showing that all lead to contradiction.

Suppose that $A$ is an exchange Mobius band. As $M$ is orientable, $A$
is one-sided.  Let $X$ be a closed regular neighborhood of $A_N$,
taken in $N$. Let $B = X \cap \bdy N$. Then $B$ is an annulus on the
boundary of the solid torus $X$.  The {\em double} of $A$, $\cover{A}$,
is the closure of $(\bdy X) \setminus B$.  Note that $\cover{A}$ is an
embedded annulus with $(\cover{A}, \bdy \cover{A}) \subset (N, \bdy
N)$. As $F$ is weakly acylindrical $\cover{A}$ must compress or
boundary-compress in $N = M \setminus \neigh_M(F)$.

Suppose that $\cover{A}$ compresses along a disk $E$. Note that $E
\cap X = \bdy E$. Compress $\cover{A}$ along $E$ to obtain a pair of
disks $C$ and $D$.  Then $S = B \cup C \cup D$ is a two-sphere
bounding $\closure{\RRPP^3 \setminus \BB^3}$ on the side which meets
$A$. Thus $M \homeo \RRPP^3$, contradicting the fact that $F$ was
two-sided and incompressible.

Suppose that $\cover{A}$ instead boundary-compresses along a disk $E$.
Let $Y$ be the closure of $\neigh_{N \setminus X}(E)$.
Boundary-compressing $\cover{A}$ along $E$ gives a disk $D$. That is,
$D$ is the closure of $(\bdy (X \cup Y)) \setminus \bdy N$. Since $F$
is incompressible and two-sided, $\bdy D$ bounds a disk $D' \subset
\bdy N$.  Note that $D'$ meets $X \cup Y$ along $\bdy D'$ only, as
$\bdy N \cap (X \cup Y)$ is nonplanar.  Now, as $F$ is incompressible,
the two-sphere $D \cup D'$ bounds a three-ball, $Z$, on the side not
meeting $A$.  It follows that $X \cup Y \cup Z$ is a solid torus with
boundary equal to a component of $\bdy N$.  However this is
impossible, as every component of $\bdy N$ has genus equal to that of
$F$.  

Thus the sum $G + H$ has no exchange Mobius bands.  Instead, suppose
that $A$ is an exchange annulus. As $M$ is orientable $A$ is
two-sided.  The weakly acylindrical hypothesis forces $A_N$ to be
compressible or boundary-compressible in $N = M \setminus
\neigh_M(F)$.

Suppose that $A_N$ is compressible along a disk $E$.  Compress $A$
along $E$ to obtain disks $C$ and $D$ with $\bdy C \cup \bdy D = \bdy
A$.  As $F$ is incompressible $\bdy C$ bounds a disk, $C' \subset F$.
Now $\bdy C = \bdy C'$ is contained inside of a seam, say
$A_r(\gamma)$.  Thus $C''$, the closure of $C' \setminus A_r(\gamma)$,
is a disk which is a union of patches and seams of $G + H$. An
innermost disk in $C''$ must be a disk patch, contradicting
Lemma~\ref{NoDiskPatches}.

Finally, suppose that $A_N$ is boundary-compressible along a disk $E$.
Boundary-compress $A_N$ to obtain a disk $D$.  By incompressibility
$\bdy D \subset \bdy N$ bounds a disk $D' \subset \bdy N$.  Let $Z$ be
the three-ball with boundary $D \cup D'$.  If $E \subset Z$ then $A_N$
is compressible (see Figure~\ref{Tube}) yielding contradiction as in
the proceeding paragraph.

\begin{figure}
\psfrag{E}{$E$}
\psfrag{A}{$A_N$}
\psfrag{F}{$\bdy N$}
$$\begin{array}{c}
\epsfig{file=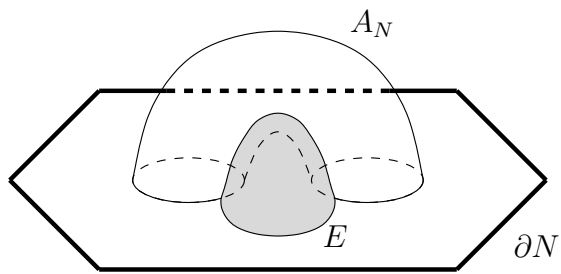, height = 3.5 cm}
\end{array}$$
\caption{A boundary-compressible tube}
\label{Tube}
\end{figure}

If $E \cap Z = \emptyset$ then the situation is more delicate.  Let
$W$ be the closure of the component of $N \setminus A_N$ containing
$E$ and $Z$.  Then $W$ is a solid torus having $E$ as a meridional
disk.  Note that $A_N \subset \bdy W$.  

Let $X$ be the closure of $\neigh_N(A_N)$.  Let $\cover{A}$ be the
closure of $\bdy X \setminus \bdy N$.  Take $F' = (\bdy N \setminus X)
\cup \cover{A}$.  Then $F' = T \cup F''$ where $T =
\bdy\left(\,\closure{W \setminus X}\,\right)$ is a torus and $F$ is
isotopic to $F''$.  (The surface $F'$ is obtained by performing an
{\em irregular exchange} along the annulus $A$.)  To check that $F''$ is
indeed isotopic to $F$ recall that $F$ is two-sided and
note that the annulus $\bdy N \cap (W \cup X)$ may be isotoped,
relative to its boundary, across the solid torus $W \cup X$.  See
Figure~\ref{ExchangeOnTunnel} for a schematic picture of the
cross-section of $W$.

\begin{figure}
\psfrag{S}{$F$}
\psfrag{S''}{$F''$}
\psfrag{A}{$A$}  
\psfrag{X}{$W \setminus X$}
$$\begin{array}{c}
\epsfig{file=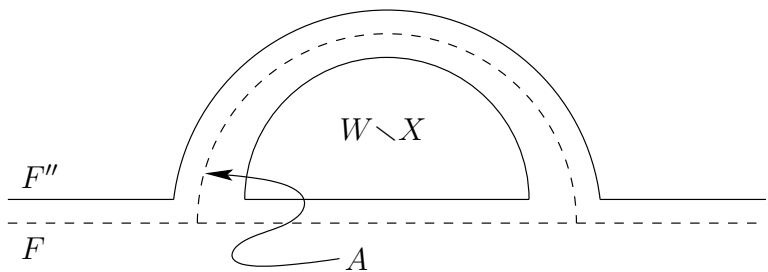, height = 3.5 cm}
\end{array}$$
\caption{An irregular exchange along $A$}
\label{ExchangeOnTunnel}
\end{figure}

Either $w(F'') < w(F)$ or $F''$ is not normal.  (Again, see
Figure~\ref{SmallNeighborhood}.) In the latter case there is an
isotopy of $F''$ to $F'''$, supported in a small neighborhood of some
face in the two-skeleton, which reduces weight by two.  See
Figure~\ref{BentArcWithIsotopy}.  Both possibilities contradict
Lemma~\ref{IncompressibleImpliesNormal} because $F$ is least weight in
its isotopy class.
\end{proof}

We now deduce our main theorem as a corollary of
Theorem~\ref{WeaklyAcylImpliesFundamental}:

\begin{theorem}
\label{A}
Suppose that $M$ is a closed, orientable three-manifold.  Then $M$
admits only finitely many inequivalent surface bundle structures with
strongly irreducible monodromy.
\end{theorem}

\begin{proof}
It follows from Theorem~\ref{WeaklyAcylImpliesFundamental} and
Lemma~\ref{FiniteFundamental} that there are finitely many isotopy
classes of weakly acylindrical surfaces in $M$.  As in
Remark~\ref{WeaklyAcyl}, any bundle with strongly irreducible
monodromy has a weakly acylindrical fibre.  Finally, a given closed
embedded surface in $M$ is isotopic to the fibre of at most two bundle
structures, up to equivalence.  This gives the desired conclusion.
\end{proof}

Now suppose that $M$ is a closed, orientable three-manifold with
triangulation $T$.  Let $\{F_i\}$ be the fundamental normal surfaces
in $T$ with negative Euler characteristic while $\{T_i\}$ and
$\{P_i\}$ are those fundamental surfaces of Euler characteristic zero
and positive, respectively.  By Lemma~\ref{FiniteFundamental} each of
these collections is finite.  Let $K = |\{F_i\}|$ be the number of
fundamental surfaces with negative Euler characteristic.  Let $P =
\max\{-\chi(F_i)\}$.

We end this section by sketching a proof of:

\begin{theorem}
\label{B}
If $(h,\phi)$ is a surface bundle structure on $M$ with genus $g > 1$
then $h^i$ is weakly reducible for all integers $i$ where $1 \leq i
\leq c(M,T) \cdot (2 g - 2)$ and $c(M,T) = 1/(3KP)$.
\end{theorem}

\begin{proof}
Suppose that $F \subset M$ is a fibre of a surface bundle structure
with monodromy $h$.  Isotope $F$ to be normal with respect to the
triangulation $T$ and the least weight such.  Suppose that $F$
decomposes as a Haken sum.  Then, by Theorem~2.2
of~\cite{JacoOertel84}, $F = \sum n_i F_i + \sum m_i T_i$ where the
summands with nonzero coefficient are fundamental, incompressible
normal surfaces.  (Note that no fundamental surface with positive
Euler characteristic appears as a summand; this follows directly from
Lemma~\ref{NoDiskPatches}.)

Recall that $g(F) = \genus(F) \geq 2$.  Thus some of the $n_i$ are
nonzero.  Reindex to obtain $n_1 \geq n_i$ for all $i$.  If $F_1$ is
two-sided then rewrite the Haken sum as $F = n G + H$, with $n = n_1$,
$G = F_1$, and $H$ equal to the sum of the remaining terms.  If $F_1$
is one-sided then take $F = n G + H$ where $n$ is the integer part of
$n_1/2$, $G = 2 F_1$ is the double of $F_1$, and $H$ is
the sum of the remaining terms.

Recall that Euler characteristic is additive under Haken sum.  Now, if
$n = 0$ then $n_i = 0$ or $1$ for all $i$.  It follows that $-\chi(F)
\leq K P$ and thus $-\chi(F) \cdot c(M,T) < 1$.  In this case the
theorem is trivially satisfied.

Assume from now on that $n$ is positive.  An easy estimate shows that
$n \geq -\chi(F) / (3 K P)$.  Make $F = nG + H$ a reduced sum by
isotoping $G$ and $H$ if necessary.  The surface $n G$ is $n$ normally
parallel copies of $G$.  Label these, in order, $G_1, G_2, \ldots,
G_n$.  The surface $H$ is nonempty and intersects $G$ because $F$ is
connected.  So $H \cap n G$ decomposes into parallel (in $H$) families
of curves, each family of size $n$.  Let $\{\gamma_i\}_{i = 1}^n$ be
one such family.  Each $\gamma_i \subset G_i$ yields an exchange
annulus for the sum $F = n G + H$.  An argument identical to the proof
of Theorem~\ref{A} shows these annuli to be essential.

Pick $m \leq n$.  The exchange annuli will lift to the $m^{\rm{th}}$
cyclic cover of the surface bundle.  Thus $h^m$ is not strongly
irreducible.
\end{proof}

\section{Questions}
\label{Questions}

Here are several questions, not all of which are necessarily difficult:

\begin{itemize}
\item 
Can a periodic surface automorphism be irreducible but not strongly
irreducible? 

\item
Are weakly acylindrical surfaces vertex surfaces? (See~\cite{JacoTollefson95}
for the definition of a vertex surface.)

\item
Give an algorithm to recognize strongly irreducible surface
automorphisms or, more generally, compute translation distance.

\item
How are the strongly irreducible bundle structures on $M$ distributed
among the fibred faces of the Thurston norm ball?

\item
Is there a meaningful stabilization theory for surface bundle
structures on $M$?

\item
Suppose that $h$ is pseudo-Anosov.  What does the translation distance
of $h$ imply about the hyperbolic geometry of the mapping torus $M_h$?

\item 
Can the function $n(g)$, as in Remark~\ref{MasurMinsky}, be given more
explicitly?  Remark~\ref{LowerBoundOnN(g)} only suggests a lower
bound.
\end{itemize}

\bibliographystyle{plain}
\bibliography{bibfile}
\end{document}